\documentclass[12pt,reqno]{amsart}
\usepackage{amsmath,amsfonts,amssymb}
\usepackage[mathscr]{eucal}
\usepackage{tikz}
\usetikzlibrary{shapes,arrows}

\usepackage{Alegreya}

\newcommand{\hed}[1]{\noindent\underline{#1}\ }

\newcommand{\htext}{\bf } 
\newcommand{\ntext}{\em } 
\newcommand{\bs}{\bigskip}
\newcommand{\ms}{\medskip}

\def\R{\mathbb{R}}
\def\RP{\mathbb{R}\mathrm{P}}
\def\AA{\mathbb{A}}
\def\EE{\mathbb{E}}
\def\calF{\mathcal{F}}
\def\calFA#1{{\calF_{\AA^{#1}}}}
\def\calFE#1{{\calF_{\EE^{#1}}}}
\def\calFR#1{{\calF_{\R^{#1}}}}
\def\calT{\mathcal{T}}
\def\scrN{\mathscr{N}}
\def\gam{\mbox{\raisebox{.45ex}{$\gamma$}}}
\def\vgam{\mbox{\raisebox{.45ex}{$\mathbf{\gamma}$}}}
\newcommand{\fbar}{\overline{f}}
\newcommand{\xbar}{\overline{x}}
\newcommand{\fg}{\mathfrak{g}}

\newcommand{\spam}{\operatorname{span}}

\def\vect#1{\mathbf{#1}}
\def\vi{\vect{i}}
\def\vj{\vect{j}}
\def\vzero{\vect{0}}
\def\vx{\vect{x}}
\def\vz{\vect{z}}

\def\vb{\vect{b}}
\def\ve{\vect{e}}

\def\bq{\boldsymbol{q}}

\begin{document}
\title{Introducing the classical method of moving frames}
\author{Thomas A. Ivey}
\address{Dept. of Mathematics, College of Charleston}
\email{iveyt@cofc.edu}

\begin{abstract}
The method of moving frames (repère mobile) was used by Elie Cartan
as a way of organizing the identification of differential invariants and solving equivalence problems.   In this expository paper, we discuss how moving frames
are used to determine differential invariants of curves and surfaces under the action 
of Euclidean, affine and conformal transformation groups.
\end{abstract}

\dedicatory{\it dedicated to Peter Olver, in honor of his 70th birthday}

\maketitle

\section*{Introduction}
This paper is based on a presentation describing the `classical' method of using moving frames to identify differential invariants of submanifolds, which I gave  at the {\it Symmetry, Invariants, and their Applications} conference.  I am grateful to the organizers for inviting me to give a talk.\footnote{I am also grateful for travel support provided by NSF grant DMS-2217293.}
In my lecture, I used three extended examples -- from Euclidean, affine, and conformal geometry -- to illustrate the method.  The latter two examples will show what happens when the group is higher-dimensional 
and when the submanifolds are higher-dimensional, respectively.

My presentation was, in turn, based on parts of a series of lectures on moving frames which I gave at a workshop at Ewha University, Seoul, in July 2012.   Since these lectures, a much more extensive and well-organized resource on moving frames has appeared  
 in the form of the textbook {\sl From Frenet to Cartan} \cite{Jbook} by Jeanne Clelland.

\section{Example: Three Ways of Looking at a Plane Curve}
We'll give a thumbnail sketch of the transition from defining differential invariants of submanifolds by normalizing
Taylor series coefficients, to computing those invariants from certain specialized lifts into frame bundles, and we will
do this in the simplest of cases, the Euclidean geometric of smooth planar curves.   We'll begin with determining
the curvature of a planar curve by using Euclidean transformations to normalize the Taylor expansion of the curve. 
(This is meant to recall treatments in classical textbooks such as Blaschke \cite{Bla}, but it is also the essence of Fels-Olver's approach to 
identifying differential invariants of submanifolds under group actions \cite{FOframes}.)  Then, we introduce the moving frame as a normalized mapping into 
the Euclidean group, in terms of which the curvature appears in the Maurer-Cartan form.   It is then easy to change perspective and 
regard this mapping as a section of a frame bundle, in terms of which the normalization appears as a way of simplifying the pullbacks of canonical and connection forms.

\subsection{Curvature via Taylor Series}
Let $\gam: \R \rightarrow \R^2$ be a curve in the Euclidean plane, acted on by group
$ASO(2)$ of rotations and translations.  (Sometimes the symbol $\EE^2$ is used to denote the plane equipped with this group action.)  
Without loss of generality, we will assume the curve a graph:
$$\vgam(x) = \begin{bmatrix} x \\ f(x) \end{bmatrix}$$

We use the group to normalize Taylor series of $\gam$ at $x=x_0$, by moving it into a `standard position'.  First, we translate it so that
when $x=x_0$ the curve passes through the origin:
$$\tilde\vgam(x) = \begin{bmatrix} x-x_0 \\ f(x)-f(x_0) \end{bmatrix}$$
Next, we apply a unique rotation so that $\tilde \vgam'(x_0)$ points along the positive $x$-axis:

$$ \tilde \vgam(x) = \begin{pmatrix} \cos \theta & \sin \theta \\ -\sin\theta & \cos \theta\end{pmatrix}\begin{bmatrix} x - x_0 \\ f(x)-f(x_0) \end{bmatrix},
\qquad \theta = \arctan f'(x_0).$$
Since the group element taking our original graph to one satisfying these two conditions at $x=x_0$ is uniquely determined,
the remaining coefficients of the Taylor expansion at $x=x_0$ are {\em geometric invariants}.

To calculate these invariants in terms of the original graph,  we rewrite this `standardized curve' as a new graph
$$\begin{bmatrix} \xbar \\ \fbar(\xbar) \end{bmatrix} =  \begin{pmatrix} \cos \theta & \sin \theta \\ -\sin\theta & \cos \theta\end{pmatrix}\begin{bmatrix} x - x_0 \\ f(x)-f(x_0) \end{bmatrix}.$$
Of course the new graph satisfies $\fbar(0)=0$, $\fbar'(0)=0$, but using implicit differentiation gives the next coefficient as
$$ \fbar''(0)= \dfrac{f''(x_0)}{(1+f'(x_0)^2)^{3/2}}.$$
This second derivative of the standardized curve gives the {\htext curvature} of original curve at $(x_0, f(x_0))$.

One of the applications of geometric invariants is helping solve the `matching problem', e.g., in this case how much information do we need to decide if two planar curves are congruent by Euclidean motions.  Is the curvature function 
$\kappa(x) = f''(x)/(1+f'(x)^2)^3/2$ enough to determine if two planar curves are congruent?  This is difficult to answer using the curvature formula derived from Taylor expansions.

\subsection{Curvature via Mapping into a Group}
In the Taylor series approach, for every $x$ we associate an in the Euclidean group $G=ASO(2)$ moving $\gam(x)$ into standard position.
Now we focus on this as a mapping from the real line into $G$ and examine its derivatives.  
To calculate these, we'll present $G$ as group of $3\times 3$ matrices 
$$g = \begin{pmatrix} 1 & 0 \\ \vb & A \end{pmatrix},\qquad \vb \in \R^2, A \in SO(2).$$
This has the advantage that the action of $G$ on $\R^2$ given by $g \cdot \vx = A \vx + \vb$ can be realized via matrix multiplication:
$$ \begin{pmatrix} 1 & 0 \\ \vb & A \end{pmatrix} \begin{bmatrix} 1 \\ \vx \end{bmatrix} = \begin{bmatrix}1 \\ A \vx + \vb \end{bmatrix}, \qquad  \vx \in \R^2.$$
Actually, we'll consider the {\em inverse transformations}, which move a graph with horizontal tangent at the origin back to $\gam(x) = (x,f(x))$, for which we have 
\begin{equation}\label{graphmovingframe}
\vb= \begin{bmatrix}x \\ f(x)\end{bmatrix}, \quad A=  \begin{pmatrix} \cos \theta & -\sin \theta \\ \sin\theta & \cos \theta\end{pmatrix},\qquad \theta(x) = \arctan f'(x).
\end{equation}
Computing the left-invariant derivative of $x \mapsto g \in G$ yields an element in the corresponding Lie algebra $\fg$:
$$g^{-1} \dfrac{dg}{dx} = \sec \theta \begin{pmatrix} 0 & 0 & 0 \\ 1 & 0 & -\kappa(x) \\ 0 & \kappa(x) & 0 \end{pmatrix}
$$
In particular, if we re-parametrize the mapping by arclength $s = \int \sec \theta \,dx$, then two curves $\gam_1$, $\gam_2$ are {\htext congruent via action of $G$ if and only if their curvature functions match} (up to a shift in $s$) as functions of arclength:
$$\kappa_1(s+c) = \kappa_2(s), \qquad c\in \R.$$
(This last result follows from a theorem of Cartan that two mappings $f, \tilde f: M \to G$ from manifold $M$ to Lie group $G$
differ by left-multiplication by a fixed element $a\in G$ if and only if the pullbacks of the Maurer-Cartan form under $f$ and $\tilde f$ are identical; see Thm. 1.6.12 in \cite{CfB2}.)

Why is this mapping $x\mapsto g(x)$ called a {\em moving frame}?
Let $\vi, \vj$ be unit vectors at the origin $\vzero$, tangent to the axes. 
Then
$$
\ve_1(x) = g(x)_* \vi, \quad \ve_2(x) = g(x)_* \vj
$$
gives a 1-parameter family of orthonormal pairs $(\ve_1, \ve_2)$, comprising a basis (or `frame') for the tangent space to $\R^2$ 
at $\gam(x)=g(x)\cdot \vzero$; as $x$ varies, this frame `moves' so that $\ve_1(x)$ always tangent to the curve.

\bs
\noindent
{\sl In this way, $x\mapsto g(x)$ can be identified with a {\em lift} of $\gam$ into the orthonormal frame bundle of $\R^2$, which we will discuss next.}

\subsection{Curvature via the Frame Bundle}
Now we turn to the third way of defining the curvature of a planar curve, which will get to the heart of our view of Cartan's approach
to moving frames.  To support this, as well as generalizations such as Cartan's method of equivalence, we first need to introduce frame bundles.
\subsubsection*{The General Linear Frame Bundle of $\R^n$}
Let $\calFR{n}$ be the set of tuples $\{ (\vb, \ve_1, \ldots, \ve_n) \}$ where $\vb$ is a point in $\R^n$
and $(\ve_1, \ldots, \ve_n)$ is any basis for $T_{\vb} \R^n$.
Then $\calFR{n}$ is a smooth manifold, with the structure of a {\em principal bundle} over $\R^n$ with basepoint map 
$$\pi:(\vb, \ve_1, \ldots, \ve_n) \mapsto \vb$$
and fiber $GL(n)$.  Regarding $\vb, \ve_i$ as $\R^n$-valued fns on $\calF$, we express their exterior derivatives
in terms of the basis $\{ \ve_i \}$:
\begin{equation}\label{defforms}
d\vb = \ve_i \otimes \omega^i, \quad d\ve_j = \ve_i \otimes \varphi^i_j,
\end{equation}
thus defining $n$ canonical 1-forms $\omega^i$ and $n^2$ connection 1-forms $\varphi^i_j$ on $\calF$.
(We will drop the tensor product sign from now on, when expressing vector-valued 1-forms.)
These $n+n^2$ are pointwise linearly independent, giving a globally-defined coframe on $\calF$, and satisfy {\em structure equations}
\begin{equation}\label{genstruc}
d\omega^i = - \phi^i_j \wedge \omega^j, \qquad d\varphi^i_j = - \varphi^i_k \wedge \varphi^k_j.
\end{equation}

The {\ntext Euclidean} frame bundle $\calFE{n}$ is the sub-bundle of $\calF$ where bases $(\ve_1, \ldots, \ve_n)$ are oriented
and orthonormal.  This is a principal bundle with structure group $SO(n)$. The same structure equations hold by restriction, but only $\binom{n}2$ of the connection forms are independent, since $\varphi^i_j = -\varphi^j_i$.
(Likewise, any Riemannian manifold has an orthonormal frame bundle with canonical and connection forms, but the second structure equation in \eqref{genstruc} also has a curvature term.)

\subsubsection*{Moving Frame as Lift into a Frame Bundle}
In the previous section, we identified the mapping $x \mapsto g(x)$ with an assignment of an orthonormal frame $(\ve_1, \ve_2)$
at the point $\gam(x)$.  Now we formalize that assignment as a {\em lift} $\Gamma: \R \rightarrow \calFE{2}$ defined by
\begin{equation}\label{myplanarframe}
\vb = \gam(x), \qquad \ve_1 = \begin{bmatrix} \cos \theta \\ \sin \theta \end{bmatrix}, \qquad \ve_2 = \begin{bmatrix} -\sin \theta \\ \cos\theta\end{bmatrix},
\end{equation}
where $\theta$ is as in \eqref{graphmovingframe}.  This mapping is a `lift' in the sense that it makes a commutative diagram
\begin{center}
\tikzstyle{line} = [draw, -latex', thick]
\begin{tikzpicture}
\node(F) {$\calFE{2}$};
\node[below of=F, node distance=15mm](plane) {$\R^2$};
\node[left of=plane, node distance=20mm](domain) {$\R$};
\path[line](F) -- node[right,yshift=-4mm,label={$\pi$}] {} (plane);
\path[line](domain) -- node[left] {$\Gamma$} (F);
\path[line](domain) -- node[below] {$\gam$} (plane);
\end{tikzpicture}
\end{center}
Specializing the frame bundle to $n=2$, we see that
$\calFE{2}$ is 3-dimensional with canonical forms $\omega^1, \omega^2$ and connection form $\varphi^2_1$ together
giving a global coframe.  These 1-forms are defined by the analogue of 
\eqref{defforms}
\begin{equation}\label{defforms2}
d\vb = \ve_1 \omega^1 + \ve_2 \omega^2, \qquad d\ve_1 = \ve_2 \varphi^2_1, \qquad d\ve_2 = \ve_1 \varphi^1_2 = -\ve_1 \varphi^2_1.
\end{equation}
Plugging in the expressions \eqref{myplanarframe} let us compute the pullbacks of the canonical and connection forms along the lift.
Specifically, our calculation of $g^{-1} dg/dx$ shows that the lift satisfies
$$\Gamma^* \omega^1 = ds, \quad \Gamma^*\omega^2 = 0, \quad \Gamma^* \varphi^2_1 = \kappa \,ds.$$
(Moreover, these equations uniquely characterize our choice of lift.)  Again, we observe that Euclidean invariants like the curvature and arclength differential appear as by-products of the moving frame. 

\bs
{\ntext In general, Cartan's method of moving frame consists of using the action of the fiber group to {\htext simplify as much as possible the
values of the pullbacks of the canonical and connection 1-forms}, with the goal of obtaining a lift into the relevant frame bundle that is uniquely characterized by values of these pullbacks.} 

\subsubsection*{Overview}
Next, we will illustrate Cartan's method of moving frames in the setting of a larger transformation group, while still focusing on planar curves.
Again, we will proceed by requiring more and more specific adaptations of the moving frame, with goal of identifying differential invariants of (sufficiently generic) curves.  At each step, the set of frames satisfying our chosen conditions will be sections of a \underline{reduction} $\calF_k$ of the general frame bundle to smaller fiber group.  As we pass to successive reductions $\calF_k$, it will be useful to 
keep track of which canonical and connection forms are zero (or are linearly depedent) on $\calF_k$, and which are semibasic (i.e., zero on vertical vectors).  For, each semibasic must be expressed in terms of the canonical form(s) on $\calF_k$, and the fiber group acts on these coefficient in this expression.  If the action on a particular coefficient is trivial, then that coefficient gives an invariant function on the base; 
if not, we use the group to normalize some of these coefficients, leading to the next reduction, to a smaller group preserving the chosen normalization. 

\section{Example: Affine invariants for Curves}
\subsection{The Affine Plane and Frame Bundle}
Let $ASL(2)$ be the following group of transformations of $\R^2$:
$$ \begin{pmatrix} 1 & 0 \\ \vb & A \end{pmatrix} \begin{bmatrix} 1 \\ \vx \end{bmatrix} = \begin{bmatrix}1 \\ A \vx + \vb \end{bmatrix}, \qquad \vx, \vb\in \R^2, A\in SL(2,\R).$$
These are transformations which preserve straight lines and the areas of triangles and rectangles.  In classical literature (e.g., \cite{Bla} Volume 2) they are referred to affine transformations, but are usually now referred to as equi-affine transformations to emphasize the preservation of the area form.  The symbol $\AA^2$ is sometimes used to denote the plane with the Kleinian geometry defined by this group action.

By acting on the standard basis at the origin, we again identify $ASL(2)$ with a sub-bundle $\calFA{2} \subset \calFR{2}$,
in a way that $\vb$ gives the basepoint map, and frame vectors $\ve_1,\ve_2$ are the first and second columns of $A$.
Because $\det A = 1$, on this sub-bundle the connection forms satisfy
$$\varphi^1_1 + \varphi^2_2 = 0$$
and the fiber is 3-dimensional, with structure group $SL(2,\R)$.

\ms
We'll try to apply Cartan's method, as described at the end of the last section, to derive (equi)affine invariants for planar curves.\footnote{The rest of this section constitutes a long-overdue solution of Exercise 1.7.3.1 in \cite{CfB2}.}
In other words, given a regular immersion $\gam:\R \to \R^2$, can we define a unique lift into $\calFA{2}$?
Once we identify such a lift, we should be able to solve the matching problem -- in other words, finding necessary and sufficient conditions for two such curves to be congruent by the $ASL(2)$ action.

\subsection{Adapted Lifts}
\subsubsection*{First Adaptation}
Given a regular curve $\gam:\R \to \R^2$, let $\calF_0 = \gam^* \calFA{2}$ (i.e., affine frames with basepoints on the image of $\gam$).

For oriented curves in the Euclidean, our moving frame was uniquely determined by requiring $\ve_1$ to point along the curve.
For affine curves, this turns out not to be enough to pick out a unique frame, so instead we pass to a sub-bundle containing all frames satisfying this condition
Accordingly, we define a sub-bundle $\calF_1 \subset \calF_0$ consisting of frames such that 
$\ve_1$ is tangent to $\gam$. 
Since on $\calF_0$ the differential of the basepoint is expanded as
$$d\vb = \ve_1 \omega^1 + \ve_2 \omega^2,$$
and the same equation holds on $\calF_1$, 
then $\omega^2$ pulls back to be zero on $\calF_1$.  Moreover, because $\gam$ is regular, then the 1-form $\omega^1$ is nonzero on $\calF_0$ (but zero on vertical vectors, i.e., vectors tangent to the fiber).

The fiber of $\calF_1$ is 2-dimensional.  Its structure group is the subgroup $G_1$ of upper triangular matrices $\left( \begin{smallmatrix} a & b\\ 0 & a^{-1} \end{smallmatrix} \right) \in  SL(2,\R)$ which acts on the fiber by
\begin{equation}\label{affsubact}
\ve_1 \mapsto a \ve_1, \quad \ve_2 \mapsto a^{-1} \ve_2 + b \ve_1, \qquad a,b \in \R, a\ne 0.
\end{equation}
(For $g\in G_1$, let $R_g$ denote this action.)
Thus, the {\em direction} of $\ve_1$ is fixed as we move along the fiber.  Since the equation $d\ve_1 = \ve_1 \varphi^1_1 + \ve_2 \varphi^2_1$
also holds on $\calF_1$, this means that $\varphi^2_1$ is {\ntext semibasic} (i.e., zero on vertical vectors) on $\calF_1$.
Hence there is a function $u$ on $\calF_1$ such that
$$\varphi^2_1 = u\, \omega^1.$$
(We could also have derived this relation by using the structure equations \eqref{genstruc} and $\omega^2=0$ to compute $0=d\omega^2 = -\varphi^2_1 \wedge \omega^1$, implying that $\varphi^2_1$ is a multiple of $\omega^1$ at each point of $\calF_1$.)

How does this coefficient $u$ vary along the fiber of $\calF_1$?  From \eqref{affsubact},
$$R_g^* \ve_1 = a \ve_1.$$
For a fixed $a,b$, differentiate both sides of this equation, and substitute in for $d\ve_1$ and $d\ve_2$ from \eqref{defforms}, to get
\begin{equation}\label{affcalc}
a \ve_1 R_g^* \varphi^1_1 + (a^{-1}\ve_2 +b\ve_1) R_g^* \varphi^2_1
=  a(\ve_1 \varphi^1_1 + \ve_2 \varphi^2_1).
\end{equation}
Equating the $\ve_2$ coefficients implies
$$R_g^* \varphi^2_1 = a^2 \varphi^2_1.$$
On the other hand, since $d\vb = \omega^1 \ve_1$ on $\calF_1$, and the basepoint $\vb$ is fixed by $R_g$, then 
$$R_g^*\omega^1 = a^{-1} \omega^1.$$
Comparing these last two equations implies that 
\begin{equation}\label{affuscale}
R_g^* u = a^3 u.
\end{equation}
Since $u$ only changes by scaling along the fiber, it is known as a {\em relative invariant}.

\hed{Remark}  Again, we can get the same information by differentiating $\varphi^2_1 = u\omega^1$.  Using the structure equations
(and $\varphi^1_1 + \varphi^2_2=0$ and $\omega^2=0$) yields
$$(du - 3u \omega^1_1) \wedge \omega^1=0$$
i.e., $d \log |u| \equiv 3\omega^1_1$ modulo semibasic terms.
Hence, $u$ scales like the cube of how $\ve_1$ scales along the fiber.

\subsubsection*{Second Adaptation}
What does our relative invariant $u$ mean geometrically?
It follows from \eqref{affuscale} that on each fiber $\calF_1$ either $u$ is identically zero or nowhere zero on the fiber. 
But
$$d\ve_1 = \ve_1 \varphi^1_1 + u \ve_2 \omega^1$$
shows that when $u\ne 0$ the direction $\ve_1$ of the tangent line is always changing as we move along the base.
Hence, points where $u=0$ on the fiber are {\em inflection points} of the curve traced out by the basepoint.

To continue adapting the frame, we have to choose whether we are in the case where $u$ is zero or non-zero.
Say we {\ntext assume} that we are in the generic case where $\gam$ is free of inflection points.  Then \eqref{affuscale} shows that there is a smooth sub-bundle $\calF_2 \subset \calF_1$ where $u=1$ identically.

The fiber of $\calF_2$ is 1-dimensional.  Its structure group is the nilpotent subgroup $G_2\subset G_1$ where $a=1$, and its action on the fiber is the restriction of \eqref{affsubact}.
On $\calF_2$ we have $R_g^* \omega^1= \omega^1$, so $\omega^1$ is a well-defined 1-form along $\gam$, which we designate as the
(equi)affine arclength differential.\footnote{Notice that this differential is well-defined {\em without} choosing an orientation for the curve. In fact, this differential is positive when the curve is bending counterclockwise and negative when the curve is bending clockwise.} 
For a generic parametrized curve $\gam(x)$ we can compute the affine arclength
$$s = \int \omega^1 = \int \det( \gam'(x), \gam''(x))^{1/3} dx.$$

\subsubsection*{Third (and last) Adaptation}
Since $\ve_1$ is fixed along the fiber of $\calF_2$, $\varphi^1_1$ is semibasic, so again there is a function $v$ on $\calF_2$ such that
$$\varphi^1_1 = v\,\omega^1.$$
How does $v$ vary along the fiber of $\calF_2$?
Our earlier calculation \eqref{affcalc}, with $a=1$, implies that 
$$R_g^* \varphi^1_1 = \varphi^1_1 - b \omega^1,$$
so that $R_g^*v = v -b$.
Thus, along each fiber of $\calF_2$ there is a unique point where $v=0$.  Since this uses up all our freedom to adjust the frame, we have the following:

\bs
\hed{Theorem} A generic curve in the equiaffine plane has a unique lift $\Gamma$ into $\calFA{2}$ such that
\begin{equation}\label{affrels}
\Gamma^*\omega^1 = ds, \quad \Gamma^*\omega^2 = 0, \quad \Gamma^* \varphi^2_1 = \omega^1, 
\quad \Gamma^* \varphi^1_1 = 0.
\end{equation}

\subsection{Affine Curvature}

The remaining 1-form on $\calFA{2}$ not specified by relations \eqref{affrels} must pull back via the unique lift to be a multiple of $ds$, and this multiple is a differential invariant.
In fact, the {\em (equi)affine curvature} function $\varkappa$ is defined by 
$$\Gamma^* \varphi^1_2 = -\varkappa \,ds. \qquad \text{(Blaschke, 1923)}$$

By an argument similar to that in the Euclidean case, two curves are congruent under $ASL(2,\R)$ if and only if their curvature functions coincide (up to a shift in $s$).
It also follows from \eqref{affrels} that, for a generic curve parametrized by equiaffine arclength, 
$$\gam'''(s) + \varkappa\, \gam'(s)=0.$$
When $\varkappa$ is a constant, this differential equation is easily solved, and we find that curves with constant $\varkappa$ are {\em conics}: parabolas, ellipses, and hyperbolas, depending
on whether $\varkappa=0$, $\varkappa>0$ or $\varkappa <0$ respectively.  (This perhaps explains Blaschke's choice of sign.)

\section{Example: Conformal Invariants of Curves and Surfaces}

In this section, we review the conformal geometry of the 3-sphere $S^3$, as defined 
by the group of {\em M\"obius transformations}.  Realizing this group action by matrices in $SO(4,1)$ leads us to conformal moving frames which are (surprisingly) vectors in $\R^5$.  Nevertheless, we can adapt lifts into the group (which we will identify as the conformal frame bundle) to obtain invariants for curves and surfaces.

\subsection{M\"obius Transformations and Conformal Frames}
\subsubsection*{M\"obius Transformations in $S^3$}
On $\R^5$, define the semidefinite quadratic form
$$\langle \vx, \vx \rangle = -x_0^2 + x_1^2+x_2^2+x_3^2+x_4^2.$$
The set of nonzero null vectors (i.e., satisfying $\langle \vx, \vx \rangle = 0$) form a 4-dimensional cone $\scrN \subset \R^4$.
Moreover, under the projectivization map $\pi: \R^5 \setminus \{ \vzero\} \to \RP^4$ the image $\pi(\scrN)$ is diffeomorphic to $S^3$.
(In fact, is clear that $\pi(\scrN)$ lies in the domain of the coordinate chart
$y_1=x_1/x_0$, \ldots, $y_4=x_4/x_0$ for $\RP^4$, and maps to the unit
sphere in the $y$-coordinates.)

The group of proper linear transformations preserving $\langle\, , \, \rangle$ 
is labeled as $SO(4,1)$, a 10-dimensional group.  Since it preserves the lines on $\scrN$, its action on $\R^5$ induces a well-defined action on $S^3$, which we call
the {\em M\"obius transformations}.  These are conformal (i.e., they preserve angles between tangent vectors) but they also preserve the set of circles and the set of two-dimensional spheres in $S^3$ (in fact, acting transitively on each of these).

The identification of $SO(4,1)$ with a frame bundle is easier if we choose a different presentation of the group.  We'll change to real coordinates $z_0, \ldots, z_5$ on $\R^5$ in which
\begin{equation}\label{altconform}\langle \vz, \vz \rangle = z_1^2 + z_2^2 + z_3^2 -2z_0 z_4,
\end{equation}
and let $K \simeq SO(4,1)$ denote the group of unimodular matrices preserving this form.

\subsubsection*{Conformal Frames for $S^3$}

From the \eqref{altconform} it follows that matrix $g$ belongs to $K$ if and only if  $\operatorname{det}g =1$ and its columns $(\ve_0, \ldots, \ve_4)$ satisfy
\begin{equation}\label{conferel}
\begin{aligned}
\langle \ve_0, \ve_0 \rangle =  \langle \ve_4, \ve_4 \rangle &=0, \quad
& \langle \ve_0, \ve_4 \rangle &=-1, \\
\langle \ve_i, \ve_0 \rangle = \langle \ve_i, \ve_4 \rangle &=0,
& \langle \ve_i, \ve_j \rangle &= \delta_{ij}, \qquad 1 \le i,j \le 3.
\end{aligned}
\end{equation}
In particular, its first and last columns are null vectors, and this lets us define a submersion from $K$ to $S^3$, given by 
$$\bq:g \mapsto \pi(\ve_0)$$
which gives $K$ the structure of a principal bundle over $S^3$ with a 7-dimensional structure group.   We will identify $K$ as the {\htext conformal frame bundle} $\calF$ of $S^3$.  This formalism has the advantage that the frame vectors $\ve_i$ 
take value in a fixed vector space $\R^5$.  But note that, for $p\in S^3$, a choice of 
a point in the fiber $\bq^{-1}(p)$ is more than just a choice of a basis for $T_p S^3$:
it is a choice of a lift $\ve_0$ of $p$ into the null cone, lifts of the basis into 
a choice of a second null vector $\ve_4$ satisfying $\langle \ve_0, \ve_4 \rangle =-1$,
and choices of lifts of the basis into $T_{\ve_0} \R^5$ orthogonal to both of these.

As before, define 1-forms $\omega^a_b$ on $\calF$ such that
$$d\ve_a = \ve_b \omega^b_a,\qquad 0 \le a,b \le 4.$$
These correspond to the left-invariant {\em Maurer-Cartan forms} on $K$, and satisfy the Maurer-Cartan equations $d\omega^a_b = -\omega^a_c \wedge \omega^c_b$.
But they also satisfy the linear relations 
\begin{align*}
\omega^0_4 &= 0, & \omega^4_0 &= 0, & \omega^4_4&= -\omega^0_0,\\
\omega^4_i &= \omega^i_0, & \omega^i_4 &= \omega^0_i, & \omega^j_i &= -\omega^i_j, \qquad i,j=1,2,3.
\end{align*}
implied by \eqref{conferel}.  For the sake of convenience, we'll write 
$\omega^1, \omega^2, \omega^3$ respectively for the 1-forms $\omega^1_0, \omega^2_0, \omega^3_0$, which are semibasic for $\bq$.  Because $\bq$ is a submersion, these 1-forms are pointwise linearly independent.

\subsection{Conformal Frames and Invariants of Curves}
Let $\gam:\R \to S^3$ be a regular parametrized curve, and $\calF_0 = \gam^* \calF$.  (Since the base of $\calF_0$ is one-dimensional, the pointwise span of $\omega^1, \omega^2, \omega^3$ is now only one-dimensional.)  This bundle still has 7-dimensional fibers acted on faithfully and transitive by $G$, which includes the full rotation group $SO(3)$ acting on $(\ve_1, \ve_2, \ve_3)$.
Thus, there is a sub-bundle $\calF_1 \subset \calF_0$ consisting of frames such that {\em $p_* \ve_1$ is tangent to $\gam$.}
Just as the subgroup preserving $\ve_1$ has codimension two 
in $SO(3)$, $\calF_1$ has codimension two in $\calF_0$.   We next work out the relations satisfied by the Maurer-Cartan forms on $\calF_1$.

Given an arbitrary lift $\Gamma$ into $\calF_1$, the derivative
of $\gam = \ve_0 \circ \Gamma$ has components only in the direction
of $\ve_1 \circ \Gamma$.  Comparing with 
$$d\ve_0 = \ve_0 \omega^0_0 + \ve_1 \omega^1 + \ve_2 \omega^2
+ \ve_3 \omega^3,$$
shows that $\omega^2=\omega^3=0$ on $\calF_1$.  Since the direction of $\bq_* \ve_1$ is
fixed as we move along the fiber of $\calF_1$, then   
\begin{equation}\label{confdeone}
d\ve_1 = \ve_0 \omega^0_1 + \ve_2 \omega^2_1 + \ve_3 \omega^3_1 + \ve_4 \omega^1
\end{equation}
shows that $\omega^2_1, \omega^3_1$ are semibasic on $\calF_1$, and thus are equal to some coefficients times $\omega^1$.
Given a lift $\Gamma$, we can adjust the lift pointwise by adding 
multiples of $\ve_2, \ve_3$ to $\ve_4$, so as to absorb these coefficients into the last
term in \eqref{confdeone}.  Thus, there is a codimension-two sub-bundle $\calF_2 \subset \calF_1$ where $\omega^2_1=\omega^3_1=0$.

The fibers of $\calF_2$ are 3-dimensional, and its structure group acts by
$$
\ve_0 \mapsto \lambda \ve_0, \qquad
\ve_1 \mapsto \ve_1 + \mu \ve_0,\qquad 
\ve_4 \mapsto \lambda^{-1}( \ve_4 + \mu\ve_1+\tfrac12 \mu^2 \ve_0), \qquad 
\lambda, \mu \in \R, \lambda \ne 0,
$$
in addition to the rotations $\ve_2 \mapsto \cos\theta \,\ve_2 + \sin\theta\, \ve_3$, $\ve_3 \mapsto -\sin\theta\, \ve_2 +\cos \theta\, \ve_3$.
Since on $\calF_2$ we are no longer free to add multiples of $\ve_2$ or $\ve_3$ to $\ve_4$, the 1-forms $\omega^2_4 = \omega^0_2$ and $\omega^3_4 = \omega^0_3$ are semibasic on $\calF_2$.
Moreover, the action of the structure group shows that 
$\spam\{ \ve_0, \ve_1, \ve_4\}$ is {\ntext invariant along each fiber} of $\calF_2$.
Its projection to $S^3$ gives the `osculating plane' to the curve, and curves 
in $S^3$ for which this is constant lie in a fixed affine plane in $\R^4$, i.e., they are 
 {\ntext circles} in $S^3$. 

Now we make a genericity assumption, that nowhere along $\gam$ does the curve have higher-order contact with its osculating plane. 
Then we can use the  remaining fiber group to make $\omega^2_4$ equal to $\omega^1$ and $\omega^0 = \omega^3_4=0$, resulting in a unique moving frame.

We can use the relations among the pullbacks of the Maurer-Cartan forms along this unique lift to 
derive Frenet-type equations for the moving frame. 
Letting $\omega^1=ds$ be conformal arclength differential, the coefficients of the remaining semibasic
1-forms $\omega^0_1, \omega^3_2$ occur as invariants in the Frenet equations:
$$
\dfrac{d\ve_0}{ds} = \ve_1,\quad
\dfrac{d\ve_1}{ds} = \kappa \ve_0 + \ve_4, \quad
\dfrac{d\ve_2}{ds} = \ve_0 +\tau \ve_3, \quad
\dfrac{d\ve_3}{ds} = -\tau \ve_2, \quad
\dfrac{d\ve_4}{ds} = \kappa \ve_1 +\ve_2
$$
These invariants are suggestively labeled $\kappa$ and $\tau$ because they 
have some similar properties to the analogous Euclidean invariants.
For example, conformal curves with $\tau=0$ lie on spheres in $S^3$.  However, 
 curves with $\kappa=0$ and $\tau=0$ are not circles; in fact, they conformally equivalent to {\em logarithmic spirals} in the plane, after applying stereographic projection.

\subsection{Conformal Frames and Invariants of Surfaces}
Let $\Sigma \subset S^3$ be an embedded surface, and let $\calF_0 = \calF\vert_\Sigma$ (i.e., frames whose basepoint lies on $\Sigma$).  We will briefly describe the process of constructing adapted frames for $\Sigma$.

\subsubsection*{First adaptation} Let $\calF_1 \subset \calF_0$ be sub-bundle of frames such that $\pi_* \ve_1, \pi_* \ve_2$ span the 
tangent space of $\Sigma$.  This is a principal bundle over $\Sigma$ with 5-dimensional fibers.  By analogy with the bundle $\calF_1$ in the curve case, 
the equation 
$$d\ve_0 = \ve_0 \omega^0 + \ve_1 \omega^1 + \ve_2 \omega^2
+ \ve_3 \omega^3$$
implies that $\omega^3=0$ and $\omega^1, \omega^2$ span the semibasics on $\calF_1$. 

Differentiating $\omega^3=0$ gives 
$$0 = \omega^3_1 \wedge \omega^1 + \omega^3_2 \wedge \omega^2.$$
By the Cartan Lemma (see, e.g., Lemma A.1.9 in \cite{CfB2}), there are functions $a,b,c$ on $\calF_1$ such that 
$$
\omega^3_1 = a\omega^1 + b \omega^2,\qquad
\omega^3_2 = b\omega^1 + c \omega^2.
$$
Note that if we form the analogue of the Euclidean second fundamental form,
$$\omega^3_1 \otimes \omega^1 + \omega^3_2 \otimes \omega^2 = a (\omega^1)^2 + c(\omega^2)^2 +2b \omega^1 \omega^2,$$
it turns out that this is not invariant along a fiber of $\calF_1$, even up to scaling.
However, by computing how $a,b,c$ vary along the fiber, we see that the {\ntext third fund. form}
$$\calT = b (\omega^1)^2 - b (\omega^2)^2 + (c-a) \omega^1 \omega^2$$
is well-defined (up to scaling) on fibers of $\calF_1$.
%
Considerations on $\calF_1$ lead to a few objects on $\Sigma$ which are conformally invariant:
\begin{itemize}
\item lines of curvature (which are null directions for $\calT$) and umbilic points (where $\calT$ vanishes) are conformally invariant;
\item the 2-form $\Omega = (b^2 + \tfrac14(a-c)^2) \omega^1 \wedge \omega^2$ is invariant along the fibers of $\calF_1$,
so gives a conformally invariant on $\Sigma$.  Its integral over $\Sigma$ is the famous {\em Willmore functional} ${\mathcal W}$, and surfaces that are critical for this functional are known
 as {\em Willmore immersions}.
\end{itemize}
To derive further invariants, we need to specialize the frames a bit more.

\subsubsection*{Second and Third Adaptations:} Let $\calF_2 \subset \calF_1$ be where the sub-bundle where $a+c=0$.  Then $\omega^0_3$ is semibasic
on $\calF_2$, and $\ve_3$ is fixed along the 4-d fibers.  It thus gives
the value of a well-defined mapping $\Gamma: \Sigma \to Q$, where $Q^4 \subset \R^5$ is the quadric defined by $\langle \vz, \vz \rangle = 1$.  This is known as the
{\em conformal Gauss map} and for compact surfaces $\Sigma$, 
the area of the Gauss image equals ${\mathcal W}(\Sigma)$.

Finally, we make the genericity assumption that $\Sigma$ is free of umbilic points. Then there is a {\em unique} conformal moving
frame such that $b=0$, $\omega^0_3=0$ and $\Omega = \omega^1 \wedge \omega^2$.  In terms of this frame, the vector $\ve_4$ takes value in the null cone $\scrN$, and composing with projectivizing gives a map to a well-defined {\em dual surface} in $S^3$.  Again, this has special properties with respect to the Willmore functional:

\ms
\hed{Theorem} (Bryant \cite{B84}) If $X:M^2\to S^3$ is a Willmore immersion, with non-umbilic points forming an open dense subset $U\subset M$, then its dual
immersion $\widehat{X}: U \to S^3$ is also a Willmore immersion.

\bs
\hrulefill 

\bs
Other choices of adapted frame can be useful for investigations of conformal surface geometry.  For example, consider the

\ms
\hed{Problem} Classify surfaces $\Sigma \subset \R^3$ foliated by two orthogonal families of circles.

\ms
Since the conditions are conformally invariant, it is natural to assume that 
$\Sigma$ is a surface in $S^3$ and apply conformal moving frames.
By adapting sections of $\calF_2$ to point along these circles, one can prove (see \cite{Idupin} for details):

\ms
\hed{Theorem}(Ivey) The circles must be lines of curvature (hence $\Sigma$ is a {\em cyclide of Dupin}) unless tangents to the circles
are {\ntext bisected} by the lines of curvature, in which case they are also circles, and $\Sigma$ is conformally equivalent to a Clifford torus.

\end{document}